\documentclass[a4paper,10pt]{amsart}

\numberwithin{equation}{section}
\usepackage{amsmath}
\usepackage{amssymb}
\usepackage{amsthm}
\usepackage{verbatim}
\usepackage[all]{xy}
\usepackage{enumerate}

\newtheorem{thm}{Theorem}[section]
\newtheorem{prop}[thm]{Proposition}
\newtheorem{cor}[thm]{Corollary}
\newtheorem{lem}[thm]{Lemma}

\theoremstyle{definition}

\newtheorem{dfn}[thm]{Definition}

\newtheorem{rmk}[thm]{Remark}

\newcommand{\N}{\mathcal{N}}
\newcommand{\M}{\mathcal{M}}
\newcommand{\R}{\mathcal{R}}
\newcommand{\T}{\mathbb{T}}
\newcommand{\D}{\mathbb{D}}
\newcommand{\E}{\mathbb{E}}
\newcommand{\F}{\mathbb{F}}

\newcommand{\X}{\mathcal{X}}

\newcommand{\A}{\mathcal{A}}
\newcommand{\U}{\mathcal{U}}

\newcommand{\PCEV}{\mathbb{P}}

\newcommand{\linspan}{\textrm{span}}

\newcommand{\ep}{\epsilon}
\newcommand{\Tr}{{\rm Tr}}

\newcommand{\Lp}{\mathcal{L}^p}
\newcommand{\sgn}{\textrm{sign}}

\title{The Walsh basis in the  $L^p$-spaces of hyperfinite III$_\lambda$ factors, $0 < \lambda \leq 1$}
\author{M. Caspers, D. Potapov, F. Sukochev   }

\address{M. Caspers, Radboud Universiteit Nijmegen, IMAPP, FNWI, Heyendaalseweg 135, 6525 AJ Nijmegen,
the Netherlands}
\email{caspers@math.ru.nl}
\address{D. Potapov, F. Sukochev, School of Mathematics and Statistics, UNSW, Kensington 2052, NSW, Australia}
\email{d.potapov@unsw.edu.au}
 \email{f.sukochev@unsw.edu.au}

\begin{document}

\begin{abstract} We introduce a non-commutative Walsh system and prove that it forms a Schauder basis in the $L^p$-spaces ($1 < p < \infty$) associated with the   hyperfinite III$_\lambda$-factors $(0 < \lambda \leq 1)$. 
\end{abstract}
 
\maketitle

 \thispagestyle{empty}

\section{Introduction}
 
In the present paper we study the non-commutative $L^p$-spaces associated with the hyperfinite factors of type III$_\lambda$, where $0 < \lambda \leq 1$. In particular, we are interested in the decomposition of this space in terms of $L^p$-spaces of matrix algebras and the construction of a very classical Schauder basis, namely the {\it Walsh system}.

\vspace{0.3cm}
 
Recall that the classical Walsh system is defined as follows. One firstly defines the Rademacher functions:
\[
r_j(x) = \sgn \left (\sin\left(  2^j \pi   x  \right)\right), \qquad j \in \mathbb{N}, x \in [0,1].
\] 
The {\it classical Walsh system}, see e.g. \cite{KashinSaakyan}, is defined as the sequence of functions given by: 
\begin{equation}\label{EqnWalshFunction}
w_n = \prod_{\gamma_i \not= 0} r_i, \qquad \textrm{ where } n = \sum_{i = 0}^\infty \gamma_i 2^i, \gamma_i \in \{0,1\}.
\end{equation}
It is a classical result that the sequence $(w_n)_{n \in \mathbb{N}}$ forms a Schauder basis in the spaces  $L^p([0,1], \mu)$ for every $ 1 < p < \infty$, see \cite[Theorem IV.15]{KashinSaakyan}. Here $\mu$ denotes the Lebesgue measure. 
 
Proper non-commutative generalizations of the Walsh system have been found for the $L^p$-spaces associated with the hyperfinite II$_1$ and II$_\infty$ factor \cite{SukFer}. Also, related problems have been studied in \cite{DFPS}, \cite{SukFerSymmetric}, where non-commutative trigonometric systems and non-commutative Vilenkin systems  where constructed. Furthermore, in \cite{PotSuk} a non-commutative Haar system was built for hyperfinite type III$_\lambda$ factors, $0 < \lambda \leq 1$.  

Here, we continue this line by constructing a Walsh system for the hyperfinite III$_\lambda$-factors, where $0 < \lambda \leq 1$.  We elaborate on the  special commutative case  $L^p([0,1], \mu_\alpha)$. Here, $\mu_\alpha$ is the Lebesgue measure in case $\alpha = \frac{1}{2}$. In case $\alpha \not = \frac{1}{2}$, the measure $\mu_\alpha$ is a   biased measure  which is singular to the Lebesgue measure and appears naturally in the construction of III$_\lambda$ factors, c.f. \cite{Kakutani}.

\vspace{0.3cm}

The structure of the paper is as follows. Section \ref{Preliminaries} recalls the necessary results on general non-commutative $L^p$-spaces. In Section \ref{SectHyperfinite} we introduce the hyperfinite III$_\lambda$ factors and fix notation. Section \ref{SectWalshBasis} contains our main result, which is the construction of a non-commutative Walsh system as a Schauder basis in the $L^p$-spaces associated with the hyperfinite III$_\lambda$ factors, $1 < p < \infty, 0 < \lambda < 1$. In Section \ref{SectIII1} we construct a Walsh system for the hyperfinite III$_1$ factor. Finally, we make remarks on the classical case in Section \ref{SectClassical}.

\section{Preliminaries on non-commutative $L^p$-spaces}\label{Preliminaries}

Let $\M$ be a von Neumann algebra with predual $\M_\ast$. For $\omega \in \M_\ast, x \in \M$, we write $x\omega \in M_\ast$ for the functional given by $(x\omega)(y) = \omega(yx), y \in \M$. Similarly, $\omega x \in \M_\ast$ denotes the functional $(\omega x)(y) = \omega(xy), y \in \M$.

\subsection{Non-commutative $L^p$-spaces}\label{SectLpSpaces}

Non-commutative $L^p$-spaces  appear in different guises. Haagerup \cite{HaagerupLp} as well as Connes and Hilsum \cite{Hilsum} gave different, but equivalent definitions of $L^p$-spaces associated with an arbitary von Neumann algebra. In \cite{Kosaki} Kosaki showed that for a von Neumann algebra with a faithful, normal state, the $L^p$-spaces are isometrically isomorphic to complex interpolation spaces between a von Neumann algebra and its predual. This is the point of view that is most suitable for our purposes. We recall the necessary definitions and notation here.

For the details on the complex interpolation method, we refer to \cite{BerghLof}. Let $\M$ be a von Neumann algebra with faithful, normal state $\omega$. 
We consider the non-dotted part of the (commutative) diagram:
\begin{equation}  
 \xymatrix{
 & \M_\ast \ar@{^{(}->}[dr]^{ \omega \mapsto \omega } & \\
   \M \ar@{^{(}->}[ur]^{x \mapsto x\omega}\ar@{^{(}->}[dr]_{x \mapsto x}\ar@{^{(}-->}[r]^{i^p  } & (\M, \M_\ast)_{[\frac{1}{p}]} \ar@{^{(}-->}[r] & \M_\ast.  \\
& \M \ar@{^{(}->}[ur]_{ x \mapsto x \omega}&}    \label{EqnLpInterpolation}
\end{equation}
 This turns the pair $(\M , \M_\ast)$ into a compatible couple of Banach spaces \cite[Section 2.3]{BerghLof}. The complex interpolation method at parameter $\frac{1}{p}$  gives by definition a Banach space $(\M, \M_\ast)_{[\frac{1}{p}]}$ which is a subset of $\M_\ast$, where the inclusion is a norm-decreasing map. Moreover, the complex interpolation method gives an embedding:
\[
i^p\:\!: \M \rightarrow (\M, \M_\ast)_{[\frac{1}{p}]}.
\] 
See also the dotted part of (\ref{EqnLpInterpolation}). 
It is proved in \cite{Kosaki} that the Banach space $(\M, \M_\ast)_{[\frac{1}{p}]}$  is isometrically isomorphic to the non-commutative $L^p$-spaces associated with $\M$ as were defined by Haagerup and Connnes/Hilsum. In particular, the construction is up to an isometric isomorphism independent of the choice of $\omega$. We simply set $\Lp(\M) = (\M, \M_\ast)_{[\frac{1}{p}]}$ as the non-commutative $L^p$-space associated with $\M$. The norm on $\Lp(\M)$ will be denoted by $\Vert \cdot \Vert_p$.

\begin{rmk}\label{RmkL1Equality}
We have an {\it equality} of Banach spaces $\mathcal{L}^1(\M) = \M_\ast$, see \cite[Theorem 4.2.2]{BerghLof}.  By the same argument $\M$ is isometrically isomorpic to $\mathcal{L}^\infty(\M)$ via the embedding $i^\infty$. 
\end{rmk}

\begin{rmk}
The  $L^p$-spaces we defined are also called $L^p$-spaces with respect to the {\it left} injection.
If one changes both the embeddings $\M \hookrightarrow \M_\ast$ in (\ref{EqnLpInterpolation}) by $x \mapsto \omega x$, the interpolated spaces are isometrically isomorphic to the present $L^p$-spaces and we refer to this construction as $L^p$-spaces with respect to the {\it right} injection. 
Other injections have been given in \cite{Kosaki}. However, the constructions in the present paper only work for the left injection and in slightly different form also for the right injection. We comment on this when it feels appropriate. Unless stated otherwise, every $L^p$-space should be understood with respect to the left injection.  
\end{rmk}

Suppose that $\N$ is a von Neumann subalgebra of $\M$ such that there exists a $\omega$-preserving conditional expectation value $\E: \M \rightarrow \N$ \cite[Definition IX.4.1]{TakII}. Denote the inclusion by $j: \N \rightarrow \M$. 
 Let $\E': \M_\ast \rightarrow \N_\ast: \omega \mapsto \omega\vert_\N$ be the restriction map and also consider the extension map $j': \N_\ast \rightarrow  \M_\ast: \omega \mapsto \omega \circ \E$.   Note that,
\begin{eqnarray}
\qquad (\E(x) \omega)(y) =  \omega(y\E(x)) = \omega(yx) = (x \omega)(y), &\qquad& x \in \M, y \in \N, \label{EqnCompatibleI}\\
(x \omega)(y)  = \omega(yx) = \omega(\E(y)x) = (x\omega)(\E(y)), &\qquad& x \in \N, y \in \M. \label{EqnCompatibleII} 
\end{eqnarray}
It follows from (\ref{EqnCompatibleI}) that the pair given by $\E$ and $\E '$ forms a morphism in the category of compatible couples of Banach spaces \cite[Section 2.3]{BerghLof} (which means by definition that $i^1_{\N} \circ \E = \E' \circ i^1_{\M}$, where $i^1_{\N}$ and $i^1_{\M}$ denote the map $i^1$ of (\ref{EqnLpInterpolation}) for respectively $\N$ and $\M$, see also Remark \ref{RmkL1Equality}).
By complex interpolation, we obtain a norm-decreasing map:
\begin{equation}\label{EqnInterpolatedCEV}
\E^p: \Lp(\M) \rightarrow \Lp(\N), \qquad 1 \leq p \leq \infty.
\end{equation} 
It follows from (\ref{EqnCompatibleII}) that the pair given by $j$ and $j'$ forms a morphism in the category of compatible couples of Banach spaces \cite[Section 2.3]{BerghLof}. Complex interpolation yields a norm-decreasing map:
\[
j^p: \Lp(\N) \rightarrow \Lp(\M).
\]
In fact, $j^p$ is isometric, since 
\[
\Vert x \Vert_p =  \Vert \E^p \circ j^p(x)\Vert_p \leq   \Vert j^p(x) \Vert_p \leq \Vert x \Vert_p, \qquad x \in \Lp(\N).
\]
Hence, we may {\it identify} $\Lp(\N)$ as a 1-complemented closed subspace of $\Lp(\M)$. 

\begin{rmk} \label{RmkLeftMultiplication}
Also {\it left} multiplication is compatible with respect to (\ref{EqnLpInterpolation}), i.e. 
\[
 x (y\omega) = (xy) \omega, \qquad x,y \in \M.
\]
Therefore, for every $x \in \M$, we can interpolate left multiplication with $x$ to give a bounded map $m_x^p$ determined by
\begin{equation}\label{EqnLeftMultiplication}
m_x^p: \Lp(\M) \rightarrow \Lp(\M): i^p(y) \mapsto i^p(xy), \qquad y \in \M.
\end{equation}
For $x \in \M$ and $y \in \Lp(\M)$, we conveniently write $xy$ for $m_x^p(y)$. 
\end{rmk}

\subsection{Martingales}\label{SectMartingales}
Let $\M$ be a von Neumann algebra with faithful, normal state $\omega$ as in Section \ref{SectLpSpaces}. Let $(\M_s)_{s\in \mathbb{N}}$ be an increasing filtration of von Neumann subalgebras of $\M$ such that their union is $\sigma$-weakly dense in $\M$. Suppose that there exist $\omega$-preserving conditional expectation values $\E_s: \M \rightarrow \M_s$. Define $\D_s = \E_s - \E_{s-1}$. As explained, we get a sequence of 1-complemented closed subspaces of $\Lp(\M)$,
\[
\Lp(\M_0) \subseteq \Lp(\M_1) \subseteq \Lp(\M_2) \subseteq \ldots \subseteq \Lp(\M),
\]
with projections $\E_s^p: \Lp(\M) \rightarrow \Lp(\M_s)$ and differences $\D_s^p = \E_s^p - \E_{s-1}^p$.

A {\it $L^p$-martingale} with respect to $(\M_s)_{s\in \mathbb{N}}$ is a  sequence $(x_s)_{s \in \mathbb{N}}$ with $x_s \in \Lp(\M)$ and $\E_s^p(x_{s+1}) = x_s$. In particular $x_s \in \Lp(\M_s)$ and $x_{s} - x_{s-1} = \D_s^p(x_s)$. 
A $L^p$-maringale $(x_s)_{s \in \mathbb{N}}$  is {\it finite} if there is a $n \in \mathbb{N}$ such that for all $s \geq n$ we have $\D_s^p(x_s) = 0$. If $x \in \Lp(\M)$, then the sequence $(\E_s^p(x))_{s\in \mathbb{N}}$ is a $L^p$-martingale. Such sequences are called {\it bounded} $L^p$-martingales. Note that the original definition of bounded is different, see \cite[Remark 6.1]{HaagerupJungeXu}. It follows that finite $L^p$-martingales are bounded. 

The following theorem follows from the Burkholder-Gundy inequalities, as first proved in the present setting in \cite{JungeXu}. The theorem also appears in \cite{HaagerupJungeXu}, where the notation is closer to ours.

\begin{thm}[Theorem 6.3 of \cite{HaagerupJungeXu}]\label{ThmUnconditionalMartingales}
Let $1 < p < \infty$. There exists a constant $C_p$, such that for every finite $L^p$-martingale $(x_s)_{s \in \mathbb{N}}$ and every choice of signs $\ep_s \in \{ -1, 1\}$, 
\[
\Vert \sum_{s =0}^\infty \ep_s \D_s^p(x_s) \Vert_p \leq C_p 
\Vert \sum_{s =0}^\infty \D_s^p(x_s) \Vert_p.
\] 
\end{thm}
It follows directly that the statement holds for every bounded $L^p$-martingale. 

\section{The setup: non-commutative $L^p$-spaces associated with the hyperfinite factors} \label{SectHyperfinite}
In this section, we fix the notation for the rest of this paper. We introduce hyperfinite factors as the direct limit of matrix algebras. 

\subsection{Hyperfinite factors}
The results in this section can be found in \cite{TakII} and \cite{TakIII}. Let $\mathbb{N}$ denote the natural numbers including 0. 
We define the following matrix algebras:
\[
\N_s =
\left\{ 
\begin{array}{ll}
 M_2(\mathbb{C})^{\otimes \frac{s+1}{2}}&\textrm{ for } s \in 2\mathbb{N}+1, \\
  M_2(\mathbb{C})^{\otimes \frac{s}{2}} \bigotimes
\left(
\begin{array}{cc}
\mathbb{C} & 0 \\
0 & \mathbb{C}
\end{array}
\right)&\textrm{ for } s \in 2\mathbb{N}.
\end{array}
\right.
\]
For $s \in 2 \mathbb{N}+1$, we consider $\N_s$ as a subalgebra of $\N_{s+1}$ by means of the embedding $x \mapsto x \otimes 1$. For $s \in 2\mathbb{N}$, there is a natural inclusion $\N_s \subseteq \N_{s+1}$.  Fix $0 < \alpha \leq \frac{1}{2}$ and let:
\[
A_1 = \left(
\begin{array}{cc}
\alpha & 0 \\
0 & 1- \alpha
\end{array}
\right), \qquad A_n = A_1^{\otimes n}, n \in \mathbb{N}.
\]
Define a state $\rho_s$ on $\N_s$ by setting 
\[
\rho_s( x ) = \Tr(x A_{\lceil \frac{s+1}{2} \rceil}), \qquad x \in \N_s.
\]
Here $\lceil \frac{s+1}{2} \rceil$ is the smallest integer that is greater than or equal to $\frac{s+1}{2}$.  

We let $\R_\alpha$ be the von Neumann algebra given by the infinite tensor product of $M_2(\mathbb{C})$ equipped with the states $\rho_1$, see \cite[Section XVIII.1]{TakIII}. 
Then, $\R_\alpha$ is a type III$_\lambda$ factor where $\lambda = \frac{\alpha}{1-\alpha}$ in case $0 < \alpha < \frac{1}{2}$ and  $\R_\alpha$ is factor of type II$_1$ in case $\alpha = \frac{1}{2}$. We have natural injective $\ast$-homomorphisms 
\[
\pi_s: \N_s \rightarrow \R_\alpha, \qquad s \in \mathbb{N}.
\]
 Furthermore, there is a distinguished faithful normal state $\rho_\alpha$ on $\R_\alpha$, which is characterized by the property:
\begin{equation}\label{EqnRhoCharacterization}
\rho_\alpha( \pi_s(x)) = \rho_s(x), \qquad s \in \mathbb{N}, \: x \in \N_s.
\end{equation}
Moreover, we have the following lemma, which is well known.
\begin{lem}\label{LemSubfactorProperties} For every $s \in \mathbb{N}$ the following holds.
\begin{enumerate}
\item The embedding $\N_s \rightarrow \N_{s+ 1}$ carries to the inclusion $\pi_s(\N_s) \subseteq \pi_{s+ 1}(\N_{s+1})$. 
\item\label{ItemInvariance} The modular automorphism group $\sigma^{\rho_\alpha}$ leaves $\pi_s(\N_s)$ globally invariant, i.e. $\sigma^{\rho_\alpha}(\pi_s(\N_s))=\pi_s(\N_s)$. 
\item The union $\cup_{s \in  \mathbb{N}} \pi_s(\N_s)$ is $\sigma$-weakly dense in $\R_\alpha$.  
\end{enumerate}
\end{lem}

For convenience of notation, we will {\it identify} $\N_s$ with its image under $\pi_s$, so that $\N_s$ is a von Neumann subalgebra of $\R_\alpha$. By (\ref{EqnRhoCharacterization}) we see that $\rho_s$ is the restriction of $\rho_\alpha$ to $\N_s$. Property (\ref{ItemInvariance}) of Lemma \ref{LemSubfactorProperties} implies that there is a $\rho_\alpha$-preserving conditional expectation value, c.f. \cite[Theorem IX.4.2]{TakII}. From now on, we use the following notation for this map:
\begin{equation}\label{EqnEs}
\E_s: \R_\alpha \rightarrow \N_s.
\end{equation}
In addition, we set $\N_{-1} = \mathbb{C} 1$,  the one-dimensional subalgebra generated by the unit of $\R_\alpha$. We set $\E_{-1}: \R_\alpha \rightarrow \N_{-1}$ as the corresponding $\rho_\alpha$-preserving conditional expectation value, which in fact is given by the map $\rho_\alpha$.

\subsection{$L^p$-spaces associated with hyperfinite factors}
It follows from the preliminaries in Section \ref{SectLpSpaces} that we get non-commutative $L^p$-spaces $\Lp(\N_s)$, with respect to the faithful, normal state $\rho_s$, with $s \in \mathbb{N}$. Similarly, we will use the notation $\Lp(\R_\alpha)$ for the $L^p$-space associated with $\R_\alpha$ with respect to $\rho_\alpha$. 
As explained, we  {\it identify} $\Lp(\N_s)$ as a closed subspace of $\Lp(\N_{s+1})$. Similarly, we may {\it identify} $\Lp(\N_s)$ as a closed subspace of $\Lp(\R_\alpha)$ and we get a chain of  closed subspaces,
\begin{equation}\label{EqnDirectLimitSequenceI}
\Lp(\N_0) \subseteq \Lp(\N_1) \subseteq \Lp(\N_2) \subseteq \ldots \subseteq \Lp(\R_\alpha), \qquad 1 \leq p \leq \infty. 
\end{equation}
As a vector space $\N_s$ is isomorphic to $\Lp(\N_s)$ by means of the mapping $i^p$, see (\ref{EqnLpInterpolation}). For $x \in \N_s$, the norm of $i^p(x) \in \Lp(\N_s)$ may be directly computed as 
\[
\Vert i^p(x) \Vert_p = {\rm Tr}( \vert x A_s^{\frac{1}{p}} \vert^p)^{\frac{1}{p}},
\]
 see \cite[Remark 3.1]{PotSuk}.

Interpolating the conditional expectation values $\E_s$, we find projections 
\[
\E_s^p: \Lp(\R_\alpha) \rightarrow \Lp(\N_s).
\]
We will need the following approximation result.  
\begin{prop}[Theorem 8 of \cite{Goldstein}]  \label{PropApproximation}
For $1 \leq p < \infty$ and  $x \in \Lp(\R_\alpha)$,
\begin{equation}\label{EqnApproximation}
\Vert x - \E^p_s(x) \Vert_p \rightarrow 0, \qquad \textrm{ as } s \rightarrow \infty.
\end{equation}
In particular, for $1 \leq p < \infty$, the union $\cup_{s \in \mathbb{N}} \Lp(\N_s)$ is dense in $\Lp(\R_\alpha)$.
\end{prop}

\section{Non-commutative Walsh system}\label{SectWalshBasis}
 
Let $\X$ be a (complex) Banach space. Recall that a sequence ${\bf x} = (x_i)_{i \in \mathbb{N}}$ in $\X$ is called a {\it Schauder basis} if for every $x \in \X$  there are unique scalars $\alpha_i \in \mathbb{C}$ such that  $ x = \sum_{i = 0}^\infty \alpha_i x_i$. In fact, ${\bf x}$ will form a Schauder basis of $\X$ if and only if the linear span of $x_i, i \in \mathbb{N}$ is dense in $\X$ and there is a constant $C$ such that for every choice of scalars $\alpha_i \in \mathbb{C}$ and every $n, m \in \mathbb{N}$ with $n > m$,
\begin{equation}\label{EqnBasisConstant}
\Vert \sum_{i = 0}^m \alpha_i x_i \Vert_{\X} \leq C \Vert \sum_{i = 0}^n \alpha_i x_i \Vert_{\X}.
\end{equation}
The constant $C$ is also called the {\it basis constant} \cite[Section 1.a]{LindenstraussTzafriri}. 

In \cite{SukFer}, a non-commutative {\it Walsh system} was given for the $L^p$-spaces associated with the hyperfinite II$_1$-factor $\R_{\frac{1}{2}}$ for $1 < p < \infty$. Recall that the system is constructed as follows. 
Consider the matrices:
\begin{equation}\label{EqnRademacher}
r^{(0,0)} =
\left(
\begin{array}{cc}
1 & 0 \\
0 & 1 
\end{array}
\right), 
r^{(1,0)} =
\left(
\begin{array}{cc}
1 & 0 \\
0 & -1 
\end{array}
\right), 
r^{(0,1)} =
\left(
\begin{array}{cc}
0 & 1 \\
1 & 0 
\end{array}
\right), 
r^{(1,1)} =
\left(
\begin{array}{cc}
0 & 1 \\
-1 & 0 
\end{array}
\right).
\end{equation}
For $n \in \mathbb{N}$ we consider the binary decomposition $n = \sum_{i = 0}^\infty \gamma_i 2^i$, where $\gamma_i \in \{0, 1\}$. We define:
\begin{equation}\label{EqnWalshBasis}
w_n = \bigotimes_{i  =0}^\infty r^{(\gamma_{2i}, \gamma_{2i+1})}.
\end{equation}
The sequence ${\bf w} = (w_n)_{n \in \mathbb{N}}$ is called the {\it Walsh system}. For $\alpha = \frac{1}{2}$, the state $\rho_{\frac{1}{2}}$ is a trace and $\R_{\frac{1}{2}}$ is the hyperfinite II$_1$-factor. In that case, it is well-known that $\Lp(\R_{\frac{1}{2}})$ is isometrically isomorphic to the semi-finite $L^p$-spaces with respect to the trace $\rho_{\frac{1}{2}}$, see also \cite[Section 2]{HaagerupLp}. Recall that the latter space can be defined as the completion of $\R_{\frac{1}{2}}$ with respect to the norm $\Vert x \Vert_p = \rho_{\frac{1}{2}}(\vert x \vert^p)^{\frac{1}{p}}$.
\begin{thm}[Proposition 5 of \cite{SukFer}]\label{ThmWalshBasisTypeII}
For $1 < p < \infty$, the Walsh system ${\bf w}$ forms a Schauder basis in the semi-finite $L^p$-spaces associated with the trace $\rho_{\frac{1}{2}}$ on $\R_{\frac{1}{2}}$.  
\end{thm}

In the present paper, we extend the result to the hyperfinite factors $\R_\alpha$. Considering the interpolation structure as described in Section \ref{SectLpSpaces}, we can consider the Walsh system ${\bf w}$ as a sequence in $\Lp(\R_\alpha)$ by means of the embedding $i^p$, see (\ref{EqnLpInterpolation}).  
We prove that ${\bf w }$ is a Schauder basis in $\Lp(\R_{\alpha})$ for $1 < p < \infty$. 

\begin{rmk}
Note that we do not incorporate $p$ explicitly in the notation of our basis ${\bf w}$. To justify this, note that by definition $\Lp(\R_\alpha)$ is as a set a subset of $(\R_\alpha)_\ast$, though their norms are different of course. Also $\R_\alpha \simeq i^\infty (\R_\alpha)$ is identified as a subset of $(\R_\alpha)_\ast$ by means of (\ref{EqnLpInterpolation}). As an element of $(\R_\alpha)_\ast$, the definition of {\bf w} does not depend on $p$. 
The principle of this slight abuse of notation is comparable to the fact that one does not distinguish a classical Walsh function (\ref{EqnWalshFunction}) as an element of $\Lp([0,1], \mu)$ for different $p$.
\end{rmk} 

We fix some auxiliary notation. Let $s \in \mathbb{N}$.  Recall that $\E_s: \R_\alpha \rightarrow \N_s $ was defined in (\ref{EqnEs}). Put 
\[
\D_s = \E_s - \E_{s-1},
\]
 and set 
$
\U_s = \D_s(\R_\alpha)$. Note that $\U_s \subseteq \N_s$. 
Moreover,
\[
\begin{split}
 \U_s =& \:\textrm{span} \left\{ w_n \mid 2^s \leq n < 2^{s+1} \right\} \\
=&\:\left\{
\begin{array}{ll}
\linspan \{ M_2(\mathbb{C})^{\otimes \frac{s}{2}} \otimes r^{(1,0)}   \},  & \textrm{ if } s \in 2 \mathbb{N}, \\
\linspan \{ M_2(\mathbb{C})^{\otimes \frac{s-1}{2}} \otimes r^{(0,1)}, M_2(\mathbb{C})^{\otimes \frac{s-1}{2}} \otimes r^{(1,1)} \},  & \textrm{ if } s \in 2  \mathbb{N}+1.
\end{array}\right.
\end{split}
\] 
Define the {\it Rademacher} matrices:
\[
r_s = 
\left\{
\begin{array}{ll}
\left( \bigotimes_{i = 1}^{\frac{s}{2} } 1 \right) \otimes r^{(1,0)} \in \N_s, & \textrm{ if } s \in 2\mathbb{N},  \\
\left( \bigotimes_{i = 1}^{\frac{s-1}{2}} 1 \right) \otimes r^{(0,1)} \in \N_s, & \textrm{ if } s \in 2 \mathbb{N}+1.   
\end{array}
\right.
\]
In particular, $r_s \in \U_s$. 

\begin{lem}
For $n \in \mathbb{N}$ and $k \in \mathbb{N}$ such that $2^k \leq n < 2^{k+1}$, we have
\begin{equation}\label{EqnRademacherRelation} 
w_{n - 2^k} = w_n r_k = \ep \: r_k w_n.
\end{equation}
Here, $\ep \in \{ -1, 1 \}$ is positive, unless $k$ is odd and $2^k +2^{k-1} \leq n < 2^{k+1}$.
\end{lem}
\begin{proof}
Suppose that $k \in 2 \mathbb{N}$. Then,
\[
w_n =\left( \bigotimes_{i  =0}^{\frac{k}{2}-1} r^{(\gamma_{2i}, \gamma_{2i+1})}\right) \otimes r^{(1,0)} \qquad \textrm{and} \qquad r_k = \left( \bigotimes_{i  = 0}^{\frac{k}{2}-1} 1 \right) \otimes r^{(1,0)}.
\]
Hence, $w_n r_k = r_k w_n = \bigotimes_{i  =0}^{\frac{k}{2}-1} r^{(\gamma_{2i}, \gamma_{2i+1})}$. Taking into account that the binary decomposition of $n$ and $n - 2^k$ are the same except for the $k$-th digit, we see that  $w_n r_k = r_k w_n =  w_{n-2^k}$.

Now, consider the case $k \in 2\mathbb{N}+1$. If $2^k \leq n < 2^k +2^{k-1}$, then
\[
w_n =\left( \bigotimes_{i  =0}^{\frac{k-3}{2}} r^{(\gamma_{2i}, \gamma_{2i+1})}\right) \otimes r^{(0,1)} \qquad \textrm{and} \qquad r_k = \left( \bigotimes_{i  = 0}^{\frac{k-3}{2}} 1 \right) \otimes r^{(0,1)}.
\]
It follows again that  $w_n r_k = r_k w_n =  w_{n-2^k}$.  If $ 2^k +2^{k-1} \leq n <2^{k+1}$, then
\[
w_n =\left( \bigotimes_{i  =0}^{\frac{k-3}{2}} r^{(\gamma_{2i}, \gamma_{2i+1})}\right) \otimes r^{(1,1)} \qquad \textrm{and} \qquad r_k = \left( \bigotimes_{i  = 0}^{\frac{k-3}{2}} 1 \right) \otimes r^{(0,1)}.
\]
Using the fact that $  r^{(1,1)} r^{(0,1)} = -r^{(0,1)}  r^{(1,1)} = r^{(1,0)}$ we now get $w_n r_k = -r_k w_n =  w_{n-2^k}$. 
\end{proof}

Let $\PCEV_n: \cup_{i= 0}^\infty \N_i\rightarrow \R_\alpha$ be the projection determined by
\[
\PCEV_n \left( \sum_{i=0}^m \alpha_i w_i \right) = \sum_{i=0}^n \alpha_i w_i, \qquad m > n, \alpha_i \in \mathbb{C}.
\]
 Note that directly after the next proposition we extend the domain of $\PCEV_n$ to $\R_\alpha$, c.f. Remark \ref{RmkPExtension}.   

\begin{thm}\label{ThmInduction}
Fix $x = \sum_{i = 0}^m \alpha_i w_i \in \R_\alpha$ with $\alpha_i \in \mathbb{C}$. For every $n < m$:
\begin{equation}\label{EqnInductionHypotheses}
 w_n \PCEV_n(x) = \E_{-1}(w_n x) + \sum_{i \textrm{ with } \gamma_i =1} \D_{i}(w_n x),
\end{equation}
where $\gamma_i \in \{ 0,1\}$ are such that $n = \sum_{i = 0}^\infty \gamma_i 2^i$.
\end{thm}
\begin{proof}
The proof proceeds by induction to $n$. For $n = 0$, note that the summation on the right hand side of (\ref{EqnInductionHypotheses}) vanishes. We find:
\[
w_0 \PCEV_0(x) = \alpha_0 w_0 = \E_{-1}(w_0 x). 
\]

Now, suppose that (\ref{EqnInductionHypotheses}) holds for all numbers stricly smaller than $n$. Let $k$ be such that $2^k \leq n < 2^{k +1}$, so that $w_n \in \U_k$. Write $n' = n - 2^k$. Then, by (\ref{EqnRademacherRelation}) we find,
\begin{equation}\label{EqnInductionI}
w_n \PCEV_n (x) = w_n \sum_{i = 0}^{n} \alpha_i w_i = w_{n'} r_{k} \left( \sum_{i = 0}^{2^k-1} \alpha_i w_i\right)  +  w_{n'} r_k \left( \sum_{i=2^k}^n \alpha_i w_i \right).
\end{equation}
For the left summation on the right hand side, the appearance of the Rademacher $r_k$ ensures that $ w_{n'} r_{k} \left( \sum_{i = 0}^{2^k-1} \alpha_i w_i\right) \in \U_k$. Hence,
\begin{equation}\label{EqnLeftSumI}
 \begin{split} 
w_{n'} r_k \left( \sum_{i = 0}^{2^k-1} \alpha_i w_i \right) = & \D_k \left( w_{n'} r_k \left( \sum_{i = 0}^{2^k-1} \alpha_i w_i \right) \right).
\end{split}
\end{equation}
By (\ref{EqnRademacherRelation}) we have $r_k w_i \not \in \U_k$ for $2^k \leq i < m$. Thus, we can continue (\ref{EqnLeftSumI}) to get,
\begin{equation}\label{EqnLeftSum}
 \begin{split} 
w_{n'} r_k \left( \sum_{i = 0}^{2^k-1} \alpha_i w_i \right) =& \D_k \left( w_{n'} r_k \left( \sum_{i = 0}^{m} \alpha_i w_i \right) \right) = \D_k(w_n x).
\end{split}
\end{equation}

Next, consider the the right summation on the right hand side of (\ref{EqnInductionI}). Using (\ref{EqnRademacherRelation}), we find that
\[
  w_{n'} r_k  \left( \sum_{i=2^k}^{n} \alpha_i w_i \right) =
 w_{n'} \left( \sum_{i=0}^{n'} \beta_i w_i \right) =   w_{n'} \PCEV_{n'} \left( \sum_{i=0}^{n'} \beta_i w_i \right),
\]
for certain $\beta_i \in \mathbb{C}$, where in fact $\beta_{i} = \pm \alpha_{i + 2^k}$, with the sign depending on $n$ (the precise equality is irrelevant for the rest of the proof). Since $n'< n$, we continue this equation by induction. Taking into account the binary decomposition of $n' = n - 2^k$ we find,
\begin{equation}\label{EqnInductionApplied}
\begin{split}
    w_{n'} r_k \left( \sum_{i=2^k}^{n} \alpha_i w_i \right) =& \E_{-1}\left(w_{n'}  \left( \sum_{i=0}^{n'} \beta_i w_i \right)\right) + \!\!\!\!\!
\sum_{s \textrm{ with } \gamma_s =1, s \not =  k} \!\!\!\!\!\!\!\!\!  \D_s \left(  w_{n'}  \left( \sum_{i=0}^{n'} \beta_i w_i \right) \right) \\
  = &
\E_{-1}   \left(  w_{n} \left( \sum_{i=2^k}^{n} \alpha_i w_i \right) \right)  + \!\!\!\!\!  \sum_{s \textrm{ with } \gamma_s =1, s \not =  k}\!\!\!\!\!\!\!\!\! \D_s \left(   w_{n}\left( \sum_{i=2^k}^{n} \alpha_i w_i \right)\right).
\end{split}
\end{equation}
Now, note that $\E_{-1}(w_n w_i) \not = 0$ if and only if $i = n$. Furthermore, let $i > n$ and let $i = \sum_{s=0}^\infty \epsilon_s 2^s$, with $\epsilon_s \in \{0, 1\}$. Looking back at (\ref{EqnWalshBasis}), we see that $w_n w_i \in \U_j$, where $j$ is the largest number such that $\gamma_j \not = \epsilon_j$. Morever, since $i > n$ we have in fact $\gamma_j = 0$ and $\ep_j = 1$. Hence, for  $i>n$,  we have $\sum_{s \textrm{ with } \gamma_s =1, s \not =  k} \D_s   ( w_{n}   w_i)   = 0$. Using these observations, we continue (\ref{EqnInductionApplied}),
\begin{equation}\label{EqnRightSUMI}
\begin{split}
    w_{n'} r_k \left( \sum_{i=2^k}^{n} \alpha_i w_i \right) =& 
\E_{-1}  \left(  w_{n} \left( \sum_{i=0}^{m} \alpha_i w_i \right) \right)  + \!\!\!\!\!  \sum_{s \textrm{ with } \gamma_s =1, s \not =  k}\!\!\!\!\!\!\!\!\! \D_s \left(   w_{n}\left( \sum_{i=0}^{m} \alpha_i w_i \right)\right)\\
  = &
\E_{-1}  \left(  w_{n} x \right)  + \!\!\!\!\!  \sum_{s \textrm{ with } \gamma_s =1, s \not =  k}\!\!\!\!\!\!\!\!\! \D_s \left(   w_{n} x\right). 
\end{split}
\end{equation}
 It is now clear that filling in (\ref{EqnLeftSum}) and (\ref{EqnRightSUMI}) into (\ref{EqnInductionI}) yields the induction hypotheses.
\end{proof}
\begin{rmk}\label{RmkPExtension}
In particular, it follows that for a fixed $n \in \mathbb{N}$ the map $\PCEV_n$ has a unique extension to $\R_\alpha$ which is both bounded and normal. We replace the notation $\PCEV_n$ by its normal extension 
\[
\PCEV_n: \R_\alpha \rightarrow \R_\alpha.
\] 
Note that we do not claim yet that the bound of $\PCEV_n$ is uniform in $n$. In fact, this is true as  we prove in the remainder of this section.
\end{rmk}

Recall from Remark \ref{RmkLeftMultiplication} that left multiplication of an element $x \in \R_\alpha$ on $\Lp(\R_\alpha)$ can be obtained by complex interpolation.  We can also interpolate the maps $\D_s, \E_s, \PCEV_s$ to get maps 
\[
\begin{split}
\D_s^p: \Lp(\R_\alpha) \rightarrow  \Lp(\N_s),\\
 \E_s^p: \Lp(\R_\alpha) \rightarrow  \Lp(\N_s), \\
 \PCEV_s^p: \Lp(\R_\alpha) \rightarrow  \Lp(\R_\alpha), 
\end{split}
\]
where $1 \leq p \leq \infty$.  
 Now, by functoriality of the complex interpolation method, we find the following corollary. 
\begin{cor}
Let $1 \leq p \leq \infty$. For every $x \in \Lp(\R_\alpha), n \in \mathbb{N}$:
\begin{equation}
 w_n \PCEV_n^p(x) = \E_{-1}^p(w_n x) + \sum_{s \textrm{ with } \gamma_s \not = 0} \D_s^p(w_n x),
\end{equation}
where $\gamma_s \in \{ 0,1 \}$ are such that $n = \sum_{s = 0}^\infty \gamma_s 2^s$. 
\end{cor}
 
At this point it is usefull to recall the definition of a Schauder decomposition. 

\begin{dfn}[Section 1.g of \cite{LindenstraussTzafriri}]
Let $\X$ be a Banach space and let ${\bf X} = (\X_s)_{s \in \mathbb{N}}$ be a sequence of closed subspaces of $\X$. Then, ${\bf X}$ is called a {\it Schauder decomposition} if every $x \in \X$ has a unique decomposition
\begin{equation}\label{EqnSchauderDecomposition}
x = \sum_{s = 0}^\infty x_s, \qquad \textrm{ where }  \: x_s \in \X_s.
\end{equation}
\end{dfn}

\begin{lem}[Section 1.g of \cite{LindenstraussTzafriri}]\label{LemSchauderEquivalence}
A sequence ${\bf X} = (\X_s)_{s \in \mathbb{N}}$  of closed subspaces of $\X$ is a Schauder decomposition if the linear span of $\cup_{s\in\mathbb{N}} \X_s$ is dense in $\X$ and furthermore, there is a constant $C$ such that
\[
\Vert \sum_{s = 0}^n x_s \Vert_{\X} \leq C \Vert \sum_{s = 0}^m x_s \Vert_{\X},
\]
for every $x_s \in \X_s$ and $n < m$.
\end{lem}

 We also need the notion of an unconditional Schauder basis.  Let ${\bf X} = (\X_s)_{s \in \mathbb{N}}$ be a Schauder decomposition of $\X$. For $A \subseteq \mathbb{N}$, consider the projection:
\[
\T_A: \X \rightarrow \X: x = \sum_{s = 0}^\infty  x_s \mapsto \sum_{s \in A} x_s,
\]
where, of course, we mean that $x_s \in \X_s$.  
\begin{lem}[Proposition 1.c.6 and its subsequent remarks in \cite{LindenstraussTzafriri}]\label{LemUncondtionalEquivalence}
The following are equivalent:
\begin{enumerate}
\item\label{ItemUnconditionalI} For every $A \subseteq \mathbb{N}$, the map $\T_A$ is bounded.
\item\label{ItemUnconditionalII} For every $x \in \X$ with $x = \sum_{s = 0}^\infty x_s$, where $x_s \in \X_s$ and for every choice $\ep_s \in \{ -1, 1\}, s \in \mathbb{N}$, the sum 
$
\sum_{s = 0}^\infty \ep_s x_s,
$ 
is convergent. 
\end{enumerate}
Moreover, if these conditions are satisfied, then there is a constant $C$ such that for every $A \subseteq \mathbb{N}$, we have $\Vert \T_A \Vert \leq C$. 
\end{lem}
If   $(\X_s)_{s \in \mathbb{N}}$ satisfies the equivalent conditions of Lemma \ref{LemUncondtionalEquivalence}, then this sequence is called an {\it unconditional} Schauder decomposition.

Note that $(\N_s)_{s\in \mathbb{N}}$ is an increasing filtration of von Neumann algebras such that its union is $\sigma$-weakly dense in $\R_\alpha$. Moreover, $\R_\alpha$ is equipped with the  faithful, normal state $\rho_\alpha$.    Therefore,  Theorem \ref{ThmUnconditionalMartingales} may be applied  and we see that (\ref{ItemUnconditionalII}) of Lemma \ref{LemUncondtionalEquivalence} holds for the decomposition  $(\D_s^p(\Lp(\R_\alpha)))_{s\in \mathbb{N}}$. 

\begin{prop} \label{PropUnconditionality}
Let $1 < p < \infty$. Then,  $(\D_s^p(\Lp(\R_\alpha)))_{s\in \mathbb{N}}$ is an unconditional Schauder decomposition of $\Lp(\R_\alpha)$. 
\end{prop}

We are now in a position to prove the main theorem of this section. 

\begin{thm}\label{ThmWalshBasisTypeIII}
For $1 < p < \infty$, the Walsh system ${\bf w}$ forms a Schauder basis in $\Lp(\R_\alpha)$.
\end{thm}
\begin{proof}
It follows from Proposition \ref{PropApproximation} that the linear span of the Walsh system is dense in $\Lp(\R_\alpha)$. 
We have to prove that  (\ref{EqnBasisConstant}) with $\X = \Lp(\R_\alpha)$ holds for a certain $C$. Equivalently, we must prove that the projections $\PCEV_n^p$ are uniformly bounded in $n$. Recall that by Theorem \ref{ThmInduction} for $x \in \Lp(\R_\alpha), n \in \mathbb{N}$:
\begin{equation}
\PCEV_n^p(x)  = w_n \E_{-1}^p(w_n x) + w_n \sum_{s \textrm{ with } \gamma_s \not = 0} \D_s^p(w_n x),
\end{equation}
where $\gamma_s \in \{ 0,1 \}$ are such that $n = \sum_{s = 0}^\infty \gamma_s 2^s$. Now, left multiplication with $w_n$ is an isometric map on $\Lp(\R_\alpha)$. Hence, 
\begin{equation}\label{EqnFinalEqn}
\Vert
\PCEV_n^p \Vert = \Vert  \E_{-1}^p +  \sum_{s \textrm{ with } \gamma_s \not = 0} \D_s^p \Vert \leq \Vert  \E_{-1}^p\Vert +  \Vert \sum_{s \textrm{ with } \gamma_s \not = 0} \D_s^p \Vert ,
\end{equation}
Since we assumed that $1 < p < \infty$, the decomposition $(\D_s^p(\Lp(\R_\alpha)))_{s\in \mathbb{N}}$ is unconditional. Hence, it follows from Lemma \ref{LemUncondtionalEquivalence} that the right hand side of (\ref{EqnFinalEqn})  is uniformly bounded in $n$.
\end{proof}

\begin{rmk}
We would like to emphasize that   the fact that left multiplication is compatible with the left injection forms an essential step in the proof of Theorem \ref{ThmWalshBasisTypeIII}.  If one considers $L^p$-spaces with respect to the right injection, one can prove that
for $1 \leq p \leq \infty, x \in \Lp(\R_\alpha)$ and $n \in \mathbb{N}$:
\begin{equation}
 \PCEV_n^{p, \sharp}(x) w_n = \E_{-1}^{p, \sharp}(x w_n) + \sum_{s \textrm{ with } \gamma_s \not = 0} \D_s^{p, \sharp}(x w_n),
\end{equation}
where $\gamma_s \in \{ 0,1 \}$ are such that $n = \sum_{s = 0}^\infty \gamma_s 2^s$. Here, the maps 
$\D_s^{p, \sharp}, \E_s^{p, \sharp}, \PCEV_s^{p, \sharp}$ are the interpolated maps of $\D_s, \E_s, \PCEV_s$
 with respect to the right injection. Completely analogously, one can now prove that {\bf w} forms a Schauder basis in the right $L^p$-spaces. 
\end{rmk}

\section{The Walsh basis in the hyperfinite factor of type III$_1$} \label{SectIII1}

Here, we construct a Walsh basis in the $L^p$-spaces associated with the hyperfinite factor of type III$_1$. The construction follows the line of \cite[Section 7]{PotSuk}, however the arguments are   different as they rely on Section \ref{SectWalshBasis}.

Consider  arbitrary von Neumann algebras $\N$ and $\M$ with with faithful, normal states $\phi$ and $\psi$. For the modular automorphism group of $\phi \otimes \psi$, we have
\[
\sigma^{\phi \otimes \psi}_t = \sigma^{\phi}_t \otimes \sigma^{\psi}_t, \qquad t \in \mathbb{R}.
\]
Therefore, $\N \otimes 1$ is a von Neumann subalgebra of $\N \otimes \M$ that is globally invariant under $\sigma^{\phi \otimes \psi}$. There exists a $\phi \otimes \psi$-preserving conditional expectation value $\E_\N: \N \otimes \M \rightarrow \N \otimes 1$, \cite[Theorem IX.4.2]{TakII}. Suppose that ${\bf v} = (v_j)_{j\in \mathbb{N}}$ is a sequence in $\M$ with $v_j^\ast v_j$ equal to a multiple of the identity. We define maps:
\[
\F_{\N,j}(x) = (1 \otimes v_j) \E_{\N} ( (1 \otimes v_j^\ast) x), \qquad x \in \N \otimes \M. 
\] 
Since $\F_{\N,j}$ is the composition of left multiplications and $\E_\N$, we can use the complex interpolation method to get a bounded map:
\begin{equation}\label{EqnDifferences}
\F_{\N,j}^p: \Lp(\N \otimes \M) \rightarrow \Lp(\N \otimes 1) = \Lp(\N).
\end{equation}

Similarly,  we can consider a $\phi \otimes \psi$-preserving conditional expectation value $\E_\M: \N\otimes \M \rightarrow 1 \otimes \M$. If ${\bf u} = (u_i)_{i \in \mathbb{N}}$ is a sequence in $\N$ with $u_i^\ast u_i$ equal to a multiple of the identity,  then we set:
\[
\F_{\M,i}(x) = (u_i \otimes 1) \E_{\M} ( (u_i^\ast \otimes 1) x), \qquad x \in \N \otimes \M. 
\] 
Interpolating this map, yields a map $\F_{\M,i}^p: \Lp(\N \otimes \M) \rightarrow \Lp(1 \otimes \M) = \Lp(\M)$.

The following theorem  can be proved similarly as \cite[Theorem 7.1]{PotSuk}. For completeness and convenience of the reader, we give the proof. Recall that the shell enumeration is an enumeration of $\mathbb{N} \times \mathbb{N}$, which assigns to a pair $(i,j)$ the number
\[
\varphi(i,j) = 
\left\{
\begin{array}{ll}
j^2+i & \textrm{ if } i \leq j,\\
(i+1)^2 - j -1 & \textrm{ if } i > j.
\end{array}
\right.
\]

\begin{thm}\label{ThmTensor}
Let $1 \leq p \leq \infty$. Suppose that ${\bf u} = (u_i)_{i \in \mathbb{N}}$ and ${\bf v} = (v_j)_{j \in \mathbb{N}}$ are sequences of linearly independent unitaries in $\N$ and respectively $\M$. Denote the corresponding projections by $\F_{\N,j}^p$ and $\F_{\M,i}^p$ and suppose that $(\F_{\N,j}^p(\Lp(\N \otimes \M)))_{j \in \mathbb{N}}$ and $(\F_{\M,i}^p(\Lp(\N \otimes \M)))_{i\in \mathbb{N}}$ are Schauder decompositions of $\N \otimes \M$. Then, ${\bf u} \otimes {\bf v} = (u_i \otimes v_j)_{i,j \in \mathbb{N}}$ taken in the shell enumeration
 is a Schauder basis for $\Lp(\N \otimes \M)$.
\end{thm}
\begin{proof}
Let ${\bf z} = {\bf u} \otimes {\bf v}$ and write ${\bf z} = (z_k)_{k \in \mathbb{N}}$. Let $n,m \in \mathbb{N}$ be such that $n < m$ and consider the sum $\sum_{i = 0}^m \alpha_i z_i$, where $\alpha_i \in \mathbb{C}$.  Let $l \in \mathbb{N}$ be such that $l^2 \leq n < (l+1)^2$. There are two cases: either $l^2 \leq n \leq l^2+l$ or $l^2+l < n < (l+1)^2$. We treat the first case, since the second case can be handled similarly. First, we compute:
\[
\begin{split}
&\Vert \sum_{k = 0}^{n} \alpha_k z_k \Vert_p \leq \Vert \sum_{k = 0}^{l^2-1} \alpha_k z_k \Vert_p + \Vert \sum_{k = l^2}^{n} \alpha_k z_k \Vert_p \\ = &
\Vert \sum_{0 \leq i,j < l} \alpha_{\varphi(i,j)} u_i \otimes v_j \Vert_p + \Vert \sum_{i=0}^{n-l^2} \alpha_{\varphi(i,l)}   u_i \otimes v_l \Vert_p  
\end{split}
\]
For the two terms on the right hand side, we find:
\[
\begin{split}
\Vert \sum_{0 \leq i,j < l} \alpha_{\varphi(i,j)} u_i \otimes v_j \Vert_p & = \Vert \PCEV_{\M,l-1}^{p} \PCEV_{\N,l-1}^{p} \left( \sum_{k=0}^m \alpha_k z_k \right) \Vert_p, \\
\Vert \sum_{i=0}^{n-l^2} \alpha_{\varphi(i,l)}   u_i \otimes v_l \Vert_p  & = 
\Vert \F_{\M,l}^p \PCEV_{\N,n-l^2}^p \left( \sum_{k=0}^{m} \alpha_k z_k \right)\Vert_p, 
\end{split}
\]
where $\PCEV_{\N,s}^p = \sum_{i = 0}^s \F_{\N, i}^p$ and $\PCEV_{\M,s}^p = \sum_{j = 0}^s \F_{\M, j}^p$. Since we assumed  that the sequences $(\F_{\N,j}^p(\Lp(\N \otimes \M)))_{j \in \mathbb{N}}$ and $(\F_{\M,i}^p(\Lp(\N \otimes \M)))_{i\in \mathbb{N}}$ are Schauder decompositions of $\Lp(\N \otimes \M)$, the projections $\PCEV^p_{\N,s}$ and  $\PCEV^p_{\M,s}$ are uniformly bounded in $s$, c.f. Lemma \ref{LemSchauderEquivalence}. It follows that there is a constant $C$ such that relation (\ref{EqnBasisConstant}) is holds.
\end{proof}

Choose $0< \alpha, \alpha' < \frac{1}{2}$ such that $\R_\alpha$ and $\R_{\alpha'}$ are factors of type III$_\lambda$ and III$_{\lambda'}$ with $\frac{\log{\lambda}}{\log{\lambda'}} \not \in \mathbb{Q}$ and $\lambda = \frac{\alpha}{1-\alpha}, \lambda' = \frac{\alpha'}{1-\alpha'}$. In that case, the tensor product $\R_\alpha \otimes \R_{\alpha'}$ is isomorphic to the hyperfinite factor of type III$_1$, see \cite{Connes}, \cite{HaagerupBicentralizer}. Consider the Walsh basis ${\bf w}$ in $\Lp(\R_\alpha)$ and let ${\bf w'}$ be the Walsh basis in $\Lp(\R_{\alpha'})$. Let $\F_{\alpha, j}^p (= \F_{\R_\alpha, j}^p)$  be the projection constructed in (\ref{EqnDifferences})   and similarly consider $\F_{\alpha',i}^p(= \F_{\R_{\alpha'}, j}^p)$.

\begin{prop}\label{PropSchauderTensor}
Let $1 < p < \infty$. The decomposition $(\F_{\alpha, j}^p(\Lp(\R_\alpha \otimes \R_{\alpha'})))_{j \in \mathbb{N}}$ is a Schauder decomposition of $\Lp(\R_\alpha \otimes \R_{\alpha'})$. Similarly, $(\F_{\alpha', j}^p(\Lp(\R_\alpha \otimes \R_{\alpha'})))_{j \in \mathbb{N}}$ is a Schauder decomposition of $\Lp(\R_\alpha \otimes \R_{\alpha'})$.
\end{prop}
\begin{proof}
We only proof the first statement, since the second one can be proved similarly. Set $\PCEV_{\alpha, n} = \sum_{j=0}^n \F_{\alpha, j}$ and  $\PCEV_{\alpha, n}^p = \sum_{j=0}^n \F_{\alpha, j}^p$. In view of Lemma \ref{LemSchauderEquivalence}, we must prove that $\PCEV_{\alpha, n}^p$ is uniformly bounded in $n$. 

Let $m > n$. Consider an element $x = \sum_{0 \leq i,j \leq m} \alpha_{i,j} w_i \otimes w_j'$ with $\alpha_{i,j}\in \mathbb{C}$. We find
\[
\begin{split}
&\PCEV_{\alpha, n} \left(\sum_{0 \leq i,j \leq m} \alpha_{i,j} w_i \otimes w_j'\right) =
\sum_{0 \leq i \leq m,\: 0 \leq j \leq n} \alpha_{i,j} w_i \otimes w_j' \\ = &
\left( \iota \otimes \PCEV_{n}\right) \left( \sum_{0 \leq i,j \leq m} \alpha_{i,j} w_i \otimes w_j' \right)
\end{split}
\]
In particular, the normality of $\PCEV_{\alpha, n}$ implies that   $\PCEV_{\alpha, n} = (\iota \otimes \PCEV_{n})$, where $\iota$ is the identity on $\R_\alpha$.

Note that $\R_\alpha \otimes \N_s$ is a von Neumann subalgebra of $\R_\alpha \otimes \R_{\alpha'}$ that is globally invariant under the modular automorphism group of $\rho_\alpha \otimes \rho_{\alpha '}$. Let $\E_{\alpha, s}: \R_\alpha \otimes \R_{\alpha'} \rightarrow \R_\alpha \otimes \N_s$ be the associated $\rho_\alpha \otimes \rho_{\alpha'}$-preserving conditional expectation value. Consider also the $\rho_{\alpha '}$-preserving conditional expection value $\E_s: \R_{\alpha'} \rightarrow \N_s$.  Clearly, the uniqueness of $(\rho_\alpha \otimes \rho_{\alpha'})$-preserving conditional expectations implies that: 
\[
\E_{\alpha, s} = \iota \otimes \E_{s}.
\]  
Recall that we defined $\D_s = \E_s - \E_{s-1}$. Similarly, set $\D_{\alpha,s} = \E_{\alpha,s} - \E_{\alpha, s-1}$.

Now, we obtain the following equalities from Theorem \ref{ThmInduction}. 
\[
\begin{split}
&(1 \otimes w_n') \PCEV_{\alpha, n}(x) = (\iota \otimes w_n' ) (\iota \otimes  \PCEV_n) (x) \\= &
\left(\iota \otimes  \left( \E_{-1} + \sum_{i \textrm{ with } \gamma_i =1} \D_{i} \right)\right) \left( (1 \otimes w_n') x \right)\\ =&
\left( \E_{\alpha,-1} + \sum_{i \textrm{ with } \gamma_i =1} \D_{\alpha, i} \right) \left( (1 \otimes w_n')x\right),
\end{split}
\]
where $n = \sum_{i=0}^\infty \gamma_i 2^i$ with $\gamma_i \in \{0,1\}$. Interpolating this equation, and observing that left multiplication with $(1 \otimes w_n')$ is an isometric map on $\Lp(\R_\alpha \otimes \R_{\alpha'})$, we find that:
\begin{equation}\label{EqnFinalDecomposition}
\Vert \PCEV_{\alpha, n}^p \Vert = 
\Vert  \E_{\alpha,-1}^p + \sum_{i \textrm{ with } \gamma_i =1} \D_{\alpha, i}^p \Vert \leq 
\Vert  \E_{\alpha,-1}^p\Vert  + \Vert \sum_{i \textrm{ with } \gamma_i =1} \D_{\alpha, i}^p \Vert .
\end{equation}
By remarks similar to the ones preceeding Proposition \ref{PropUnconditionality}, it follows from Theorem \ref{ThmUnconditionalMartingales} that the decomposition $(\D_{\alpha, i}^p(\Lp(\R_\alpha \otimes \R_{\alpha'}))_{i \in \mathbb{N}}$ is an unconditional Schauder decomposition of $\Lp(\R_\alpha \otimes \R_{\alpha'})$. Hence, Lemma \ref{LemUncondtionalEquivalence} implies that the right hand side of (\ref{EqnFinalDecomposition}) is uniformly bounded in $n$.  
\end{proof}

Proposition \ref{PropSchauderTensor} implies that we may apply  Theorem \ref{ThmTensor}.

\begin{thm}
Let $1 < p < \infty$. The Walsh system ${\bf w} \otimes {\bf w}' = (w_i \otimes w_j')_{i,j \in \mathbb{N}}$ taken in the shell enumeration is a Schauder basis in $\Lp(\R_\alpha \otimes \R_{\alpha'})$; the $L^p$-space associated with the hyperfinite III$_1$ factor. 
\end{thm}

\begin{rmk}
In general a tensor product of two $L^p$-spaces, each with unconditional decomposition, does not produce a   $L^p$-space where the tensor product of the given decompositions is unconditional. The simplest example is a couple of Schatten classes with row and column decompositions.
\end{rmk}

\section{Classical $L^p$-spaces}\label{SectClassical}

For $s \in \mathbb{N}$, consider the diagonal subalgebra $\A_s \subseteq \N_s$. The weak closure of $\cup_{s \in \mathbb{N}}  \A_s$ in $\R_\alpha$ forms an abelian von Neumann algebra $\A_\alpha$, which is isomorphic to $L^\infty([0,1], \mu_\alpha)$. Here, $\mu_\alpha$ is the measure determined by:
\[
\mu_\alpha \left( \left[\frac{k}{2^n}, \frac{k+1}{2^n}  \right] \right) = \prod_{i=0}^{n-1} \left[ (1- \gamma_i) \alpha + \gamma_i (1 - \alpha) \right], 
\]
where $0 \leq k < 2^n$ and $\gamma_i \in \{0,1\}$ are such that $k = \sum_{i = 0 }^{n-1} \gamma_i 2^i$, see \cite[Section 12.3]{KadRin}. In particular, $\A_{\frac{1}{2}}$ is isomorphic to $L^\infty([0,1], \mu)$, where $\mu$ is the Lebesgue measure.

The modular automorphism group $\sigma$ leaves  $\cup_{s \in \mathbb{N}}  \A_s$ and hence $\A_\alpha$ invariant. 
From Section \ref{SectLpSpaces}, it follows that $\Lp(\A_\alpha)$ is a closed subspace of $\Lp(\R_\alpha)$.  Moreover, there exists a conditional expectation value $\E_{\A_\alpha}: \R_\alpha \rightarrow \A_\alpha$.
Since $\E_{\A_\alpha}$ projects on the diagonal matrices, we find that it acts on the Walsh system ${\bf w}$ by:
\[
\E_{\A_\alpha}(w_n) = 
\left\{
\begin{array}{ll}
w_n & \textrm{ if } n = \sum_{i = 0}^\infty \gamma_i 2^i \textrm{ with } \gamma_{2i+1} = 0 \textrm{ for every } i,\\
0 & \textrm{else.} 
\end{array}
\right.
\] 
Indeed, it follows from (\ref{EqnWalshBasis}) that $w_n$ is diagonal if and only if the odd digits in the binary decomposition of $n$ vanish.
Let ${\bf z}$ be the subsequence of ${\bf w}$ of vectors in the range of the projection $\E_{\A_\alpha}$.
Clearly, it follows from Theorem \ref{ThmWalshBasisTypeIII} that ${\bf z}$ forms a Schauder basis in $\Lp(\A_\alpha)$ for $1 < p < \infty$.  Explicitly, this system is constructed as follows. Recall that we defined the Rademacher matrices in (\ref{EqnRademacher}). Set:
\[
z_n = \bigotimes_{i =0}^\infty r^{(\gamma_{i}, 0)}, \qquad n = \sum_{i = 0}^\infty \gamma_i 2^i, \:\: \gamma_i \in \{0, 1\}.
\]
Then, ${\bf z} = (z_n)_{n \in \mathbb{N}}$.
\begin{cor}
Let $1 < p < \infty$. The system ${\bf z}$ forms a Schauder basis in $\Lp(\A_\alpha)$. Under the isomorphism $\Lp(\A_{\alpha}) \simeq L^p([0,1], \mu_{\alpha})$, we obtain the classical Walsh system (\ref{EqnWalshFunction}).  
\end{cor}

\end{document}